\documentclass[11pt]{article}
\usepackage{amssymb}
\usepackage{amsmath}

\newtheorem{theo}{Theorem}[section]

\newtheorem{prop}[theo]{Proposition}
\newtheorem{lemm}[theo]{Lemma}
\newtheorem{coro}[theo]{Corollary}
\newtheorem{rema}[theo]{Remark}

\def \C{{\mathbb C}}

\def \Z{{\mathbb Z}}
\def \p{{\mathbb P}}
\def \Q{{\mathbb Q}}

\newcommand{\cqfd}
{%
\mbox{}%
\nolinebreak%
\hfill%
\rule{2mm}{2mm}%
\medbreak%
\par%
}
\newfont{\gothic}{eufb10}
\date{\empty}
\begin{document}
\title{Potential density of rational points on the variety \\
of lines of a cubic fourfold}
\author{Ekaterina Amerik, Claire Voisin} \maketitle \setcounter{section}{-1}
\section{Introduction}
Recall that for a variety $X$ defined over a field $k$,    rational points
are said to be potentially dense in $X$ if for some finite extension $k'$ of $k$,
$X(k')$ is Zariski dense in $X$.
For example, this is the case if $X$ is a rational
or a unirational variety.

If $k$ is a number field, a conjecture of Lang and Vojta, proved in dimension
1 by Faltings, but still open even in dimension two,
predicts that varieties of general type  never satisfy potential density over $k$.
On the contrary, it is generally expected that when the canonical bundle $K_X$ is
negative (i.e. $X$ is a Fano variety) or trivial, rational points are
always potentially dense in $X$ (see for example \cite{C}
for a recent and precise conjecture).

In the Fano case, there are many examples
confirming this. First of all, many Fano varieties are unirational
(in dimension $\leq 2$, all of them are even rational over a finite extension of the definition field).
Furthermore, Fano threefolds such as a three-dimensional quartic
or a ``double Veronese cone" (a hypersurface of degree $6$ in $\p(1,1,1,2,3)$)
are birational to elliptic fibrations over $\p^2$, which makes
it possible to prove  potential density (see \cite{HarT}).

If the canonical bundle is trivial, the known examples are much
less convincing. It is well-known that rational points on
abelian varieties are potentially dense. But the simply-connected
case remains largely unsolved.

A number of  results concerning potential density of rational points on
$K3$ surfaces defined either over a number field or a function field
appeared in the last years.
Bogomolov and  Tschinkel \cite{Bo-Tschi}  proved potential
density of rational points for $K3$ surfaces  admitting an elliptic pencil, or an infinite
automorphisms group. Remark that these $K3$ surfaces are rather special;
in particular,
their geometric Picard number is never equal to $1$, whereas it is equal to $1$ for a
general projective $K3$.

The case
of $K3$-surfaces defined over a (complex) function field has been solved
for certain types of families
by Hassett and Tschinkel \cite{Hatschi}. They prove in particular the
existence
of $K3$-surfaces defined over $\mathbb{C}(t)$,
 whose geometric Picard group is  equal to $\mathbb{Z}$,
and which satisfy potential density. Their method can likely be adapted to
produce examples
over $\overline{\mathbb{Q}}(t)$.

However, no example of $K3$ surface defined over a number field, and
with geometric Picard group equal to $\mathbb{Z}$ is known to satisfy potential density.

In this paper, we consider the case of higher dimensional varieties which are
as close as possible to $K3$ surfaces, namely  Fano varieties of lines
$F$ of cubic $4$-folds $X$.
These are $4$-dimensional varieties with trivial canonical bundles, and which
even possess a non degenerate holomorphic $2$-form which will play an important role
in our study. Their similarity with
$K3$ surfaces is shown by the work of Beauville and Donagi
\cite{be-do}, which
shows that these varieties are deformations of the second punctual Hilbert scheme of a
$K3$ surface. However, 
for $X$ sufficiently general, the geometric Picard group of $F$ is equal to $\mathbb{Z}$.
By a result due to Terasoma, this property holds true even for many $X$'s defined over
$\mathbb{Q}$.
We prove the following result:

\begin{theo} \label{main} For many cubic $4$-folds $X$ defined over a number field, the corre- sponding
variety $F$ (which is defined over the same number field) satisfies the po- tential density
property.
\end{theo}
Here, the term ``many'' can be made precise and exactly depends on the
structure of the  set of points in a certain  moduli space satisfying Terasoma's density
condition (see \cite{terasoma}). In lemma \ref{conditionslemm}, we give precise conditions under which
the conclusions of Theorem \ref{main} hold.

Let us describe the strategy of our proof. We will use
two geometric ingredients.
 The first one is the rational self-map $\phi:F\dashrightarrow F$
 constructed in
\cite{voisinturin} (cf section \ref{preli}); in \cite{am-cam},
this map was already used to prove the potential density
of rational points for such a variety $F$ over an uncountable field
(say, a function field of a complex curve, as in \cite{Hatschi}), assuming
its geometric  Picard group is
$\mathbb{Z}$.

The proof in \cite{am-cam} shows in fact that for a general complex point
$l\in F(\mathbb{C})$, its orbit
$\{\phi^k(l),\,k\in \mathbb{N}\}$ is Zariski dense in $F$.
In the case where $F$ is defined over
a number field $k$, their proof does not say anything about points of $F$ over
 finite extensions of $k$,
because it is based on discarding certain ``bad" countable unions of
proper subvarieties, and it seems difficult to get any more precise control
over those subsets. Thus one cannot conclude from \cite{am-cam} (and neither from our proof)
that the Zariski density of the orbit under $\phi$ of a single point
 defined over a number field holds true.

The second ingredient is  the existence of a $2$-dimensional
 family of surfaces $(\Sigma_b)_{b\in B}$ in $F$, which are birationally equivalent
to  abelian surfaces (cf \cite{voisincdm} and section \ref{preli}).
Of course, for any surface
as above defined over a number field, rational points are potentially dense in it.
Furthermore $\phi$ is defined over the same number field as $F$.  It
thus suffices
to prove that for one such surface
$\Sigma_b$, the countably many surfaces
$\phi^l(\Sigma_b)$ are Zariski dense in $F$.

This is done in two steps, and our proof shows that this density result holds once a
non-vanishing condition (see Proposition \ref{plusfort})
for a certain $l$-adic Abel-Jacobi cycle class holds for the corresponding 
$\Sigma_b$. This non-vanishing is satisfied for many
$\overline{\mathbb{Q}}$-points in the moduli
space $B$ of these surfaces by an argument which is taken from \cite{ggp} (and attributed there
to Bloch and Nori).

In section \ref{sec2}, we prove that for many $\Sigma_b$'s defined over a number field, the Zariski closure
of the union of the surfaces
$\phi^l(\Sigma_b)$ contains at least a divisor (that is, $\Sigma_b$ is not
preperiodic), and this already necessitates this non-vanishing
condition.

In the last section, we study divisors in $F$ invariant under a power of $\phi$ and containing
a Zariski dense union of surfaces birationally equivalent to abelian surfaces.
We use the restriction of the holomorphic symplectic  form of $F$ to such divisors
to classify them (cf section \ref{secclassi}). The conclusion then comes in two steps :
we exclude one case by geometric considerations in section  \ref{secclassigeo} and
we show in section
 \ref{conclusion} that the only remaining case
is excluded under the non-vanishing condition above, which is of arithmetic nature.

\section{Preliminaries \label{preli}}

Let  $X\subset \mathbb{P}^5$ be a smooth cubic fourfold,
and $F$  the variety of lines of $X$; let us denote by $P\subset F\times X$
 the universal incidence correspondence, that is,
 $p:P\rightarrow F$ is the universal family of lines parameterized by $F$ and $q:P\rightarrow X$
  is the natural
 inclusion in  $X\subset\mathbb{P}^5$ on each of these lines. Recall from \cite{be-do} that
 the Abel-Jacobi map
 $$p_*q^*:H^4(X_\mathbb{C},\mathbb{Q})\rightarrow H^2(F_\mathbb{C},\mathbb{Q})$$
 is an isomorphism of rational Hodge structures (of bidegree $(-1,-1)$),
 which sends $h^2,\,h:=c_1(\mathcal{O}_X(1))$,
 to $c_1(\mathcal{L})$, where $\mathcal{L}$ is the restriction to $F$ of the Pl\"ucker line bundle.
 (In the whole paper, we will denote by $X_\mathbb{C}$ the complex manifold associated
 to a variety defined over a subfield of $\mathbb{C}$, and
 $H^*(X_\mathbb{C},R)$ will then denote the Betti cohomology of this complex manifold
 with coefficients in  $R$.)
The number $h^{3,1}(X)$ is equal to $1$, as shown by Griffiths theory, and we thus have
$h^{2,0}(F)=1$. Let $\omega=q_*p^*\alpha$ be a generator of $H^{2,0}(F)=H^0(F,\Omega_F^2)$,
where $\alpha$ generates $H^{3,1}(X)$. Then
$\omega$ is a non-degenerate $2$-form by \cite{be-do}.

On the other hand, it was proved by Terasoma
\cite{terasoma}  that many cubic fourfolds $X$ defined over
number fields satisfy the property that the Hodge classes in $H^4(X_\mathbb{C},\mathbb{Q})$
 are rational multiples
of $h^2$. (Terasoma formulates his result in terms of the non-existence
of non zero primitive algebraic cycle classes, but this is equivalent
to the non-existence
of non zero primitive Hodge classes,
as the Hodge conjecture is known to be true for cubic $4$-folds,
see \cite{zu}, \cite{comu}.) We deduce from this together with the fact that $H^2(F_\mathbb{C},\mathbb{Z})$
has no torsion and the class
$c_1(\mathcal{L})$ is not divisible
in $H^2(F_\mathbb{C},\mathbb{Z})$, that the corresponding variety $F$, which is also defined
over a number field, satisfies $Pic\,F_\mathbb{C}=\mathbb{Z}$,
with generator
 $\mathcal{L}$. Equivalently, the geometric Picard group of $F$ is equal
 to $\mathbb{Z}\mathcal{L}$.

For each smooth hyperplane section $Y\subset X$, the variety of lines contained in $Y$ is a
smooth surface $\Sigma_Y$ of class $c:=c_2(\mathcal{E})$, where $\mathcal{E}$ is the restriction
to $F$ of the universal rank $2$ bundle on the Grassmannian. Furthermore, these surfaces
are isotropic with respect to $\omega$. Indeed, denoting by
$$p_Y:P_Y\rightarrow \Sigma_Y,\,\,q_Y:P_Y\rightarrow Y$$
the incidence variety for $Y$, we have the following equality:
$$\omega_{\mid \Sigma_Y}=p_{Y*}\circ q_Y^*(\alpha_{\mid Y})\,\,{\rm in}\,\,H^{2,0}(\Sigma_Y).$$
But by Lefschetz theorem on hyperplane sections, we have $H^{3,1}(Y)=0$, and hence
$\alpha_{\mid Y}=0$. Thus $\omega_{\mid \Sigma_Y}=0$.

Recall next the observation made in \cite{voisincdm} : consider  first the
case where $Y$ has one ordinary
singular point $y\in Y$. Then from \cite{cl-gr}, we know that
the normalization of $\Sigma_Y$ is isomorphic to
the symmetric product $C^{(2)}_x$, where $C_x$ is the curve of lines in $Y$ passing through $y$.
This curve is the complete intersection of a cubic and a quadric in $\mathbb{P}^3$.
The map $C^{(2)}_x\rightarrow \Sigma_Y$ sends a pair of lines $l_1,\,l_2$ through $y$  to the
residual line $l$ of the intersection with $Y$ of the plane $<l_1,l_2>$ generated by the two meeting
lines $l_1$ and $l_2$ (this intersection already contains the two lines
$l_1,\,l_2$). The
 inverse of this map associates to a generic line $l$ in $Y$ not passing through $y$
the pair of lines which forms the residual intersection of the plane $<l,y>$ with $Y$
(this is indeed a  conic singular at $y$).

The curve $C_x$ is smooth and has genus $4$.

Assume now that $Y$ has exactly two more double points $y'$ and $y''$, and
$y'$, $y''$, $y$ are not on the same line.  Then
the curve $C_x$ also acquires two double points, corresponding to the
lines $l',\,l''$ joining $y$ respectively
to $y'$ and $y''$ (these lines are obviously contained in $Y$).

Thus the normalization of $C_x$ has genus $2$, and its second
symmetric product is birational to an abelian surface.
\begin{rema} {\rm This abelian surface is isomorphic to the intermediate Jacobian of the desingularization
$\widetilde{Y}$
of $Y$ obtained by blowing-up the nodes. It is also isomorphic via the Abel-Jacobi map to the
group $CH_1(\widetilde{Y})_{hom}$ (cf \cite{bloch}, \cite{cl-gr}).}
\end{rema}

Imposing three double points to $Y$ imposes three conditions to the corresponding hyperplane
of $\mathbb{P}^5$. Thus we get a $2$-dimensional family
$(\Sigma_b)_{b\in B}$ of such (highly singular) surfaces in $F$, and we will
denote by $(A_b)_{b\in B}$ the corresponding family of abelian surfaces, which are birationally
 equivalent
to them. Thus $A_b$ can be defined as $Alb\,\Sigma'_b$ for some desingularization
$\Sigma'_b$ of $\Sigma_b$.

Note that if $X$ is defined over a number field, so are $B$ and the family of
abelian  surfaces
$(A_b)_{b\in B}$.

We will need the following lemma:
\begin{lemm} Any smooth genus $2$ curve appears as the normalization
of the family of lines through a singular point of some cubic threefold with three
non colinear double points.
\end{lemm}
{\bf Proof.} A generic genus $2$ curve $C$ is the normalization of a complete intersection
of a smooth quadric and a cubic surface in $\mathbb{P}^3$ admitting exactly $2$ double points.
Indeed, it suffices to choose two generic line bundles of degree $3$ on $C$,
which will provide two maps from $C$ to $\mathbb{P}^1$. Then for a generic choice,
the induced map from $C$ to $\mathbb{P}^1\times \mathbb{P}^1$ has for image a nodal curve
$C'$.
This curve is in the linear system $\mid\mathcal{O}(3,3)\mid$ on
$\mathbb{P}^1\times \mathbb{P}^1$, which equivalently means that $C'$ is a complete intersection
$Q\cap S$, where $S$ has degree $3$ in $\mathbb{P}^3$ and $Q\cong \mathbb{P}^1\times \mathbb{P}^1$
is a smooth quadric. As the arithmetic genus of $C'$ is $4$, $C'$ has $2$ double points.

Now, following \cite{cl-gr}, let $Z$ be the blow-up of $\mathbb{P}^3$ along $C'$, and consider
the morphism from $Z$ to
$\mathbb{P}^4$ given by the cubic equations vanishing on $C'$: if $q,\,t$ are given choices
of quadratic, resp. cubic generating equations for the ideal of $C'$, this morphism is given by the linear
system
$$X_0q,\ldots, X_3q,t.$$

Denoting by $Y_0,\ldots,Y_4$ the homogeneous coordinates on $\mathbb{P}^4$ used above,
the image of $Z$ under this morphism is the cubic threefold
$Y$ which has
for equation
$$t(Y_0,\ldots,Y_3)=Y_4q(Y_0,\ldots, Y_3).$$
It is immediate to verify that $Y$ has $3$ ordinary non colinear  double points, one of them
being $o:=(0,\ldots,0,1)$, the others satisfying $Y_4=0$ and corresponding to the double points of
$C'$.
The curve of lines in $Y$ though $o$ identifies
to $C'$ by projection from $o$.

\cqfd

Our arguments in the next sections   will apply to those $X$ satisfying the
conclusions
of the following lemma. These $X$'s will thus satisfy the conclusion of theorem \ref{main}.
\begin{lemm} \label{conditionslemm}For many $X$'s defined over
 over a number field,  the following two properties
are satisfied:

- The corresponding variety of lines $F$ defined over the same number field has geometric Picard group
equal to $\mathbb{Z}$.

-  For  some (and in fact many) closed points $b\in B$ defined over a number field,
 the abelian surfaces
$A_b$
have geometric N\'eron-Severi group equal to $\mathbb{Z}$.
\end{lemm}
{\bf Proof.} The proof goes exactly as in Terasoma \cite{terasoma}.
Consider the family $\mathcal{M}$ parameterizing pairs
$(A,X)$ where $A$ is a principally polarized abelian variety which is the
intermediate Jacobian of the desingularization of a hyperplane section $Y$ of $X$ with
exactly three independent
nodes (with its canonical polarization).
$\mathcal{M}$ maps via the first projection
to the moduli space of principally polarized abelian varieties
$\mathcal{M}_A$, and via the second projection, to the
open set  $U\subset \mathbb{P}(H^0(\mathbb{P}^5,\mathcal{O}_{\mathbb{P}^5}(3)))$ parameterizing
smooth cubic $4$-folds.
Our previous result shows that
both maps are dominating.
Replacing
$\mathcal{M}$ and $\mathcal{M}_A$ by a finite cover if necessary, we may assume there
are universal families
$$\mathcal{A}\rightarrow \mathcal{M}_A,\,\,\mathcal{X}\rightarrow U,\,\mathcal{A}_{\mathcal{M}}
\times_{\mathcal{M}}\mathcal{X}_{\mathcal{M}}\rightarrow\mathcal{M}.
$$
Note furthermore that all families are defined over $\mathbb{Q}$.
We fix compatible base-points $o,\,o_A,\,o_U$ for $\mathcal{M},\,\mathcal{M}_A,\,U$.

For the general polarized abelian complex surface $A$,
one has $NS(A)=\mathbb{Z}\lambda$,
which is equivalent to the condition that the only Hodge classes in
$H^2(A,\mathbb{Q})$ are multiples of the polarization
$\lambda\in H^2(A,\mathbb{Q})$.

This is  a consequence of the stronger result saying that the monodromy
group
$$Im\,\rho_A:\pi_1(\mathcal{M}_{A,\mathbb{C}},o_A)\rightarrow Aut\,
H^2(A_{o_A,\mathbb{C}},\mathbb{Q})$$
does not act in a finite way on any non zero class in $\lambda^\perp$.

As the image of $\pi_1(\mathcal{M}_\mathbb{C},o)$ in
$\pi_1(\mathcal{M}_{A,\mathbb{C}},o_A)$ has finite index, the same holds for
$$Im\,\rho'_A:\pi_1(\mathcal{M}_\mathbb{C},o)\rightarrow Aut\,
H^2(A_{o,\mathbb{C}},\mathbb{Q}).$$
Similarly, there is no non zero class in
$H^4(X_{o_U,\mathbb{C}},\mathbb{Q})_{prim}$ on which the monodromy group
$$Im\,\rho_U:\pi_1(U_\mathbb{C},o_U)\rightarrow Aut\,
H^4(X_{o_U,\mathbb{C}},\mathbb{Q})_{prim}$$
acts in a finite way, which implies  the Noether-Lefschetz
theorem, saying that for general $X$, there is no non trivial  Hodge class
in $H^4(X,\mathbb{Q})_{prim}$.
It follows that
there is no non-zero class in
$H^4(X_{o_U,\mathbb{C}},\mathbb{Q})_{prim}$ on which the monodromy group
$$Im\,\rho'_U:\pi_1(\mathcal{M}_\mathbb{C},o_U)\rightarrow Aut\,
H^4(X_{o_U,\mathbb{C}},\mathbb{Q})_{prim}$$
acts as a finite group. Of course, we have the same conclusions
for the monodromy acting on cohomology with $\mathbb{Q}_l$ coefficients:
$$\rho'_A:\pi_1(\mathcal{M}_\mathbb{C},o)\rightarrow Aut\,
H^2(A_{o,\mathbb{C}},\mathbb{Q}_l),\,\rho'_U:\pi_1(\mathcal{M}_\mathbb{C},o_U)\rightarrow Aut\,
H^4(X_{o_U,\mathbb{C}},\mathbb{Q}_l)_{prim}.$$
Finally, we use the density theorem of Terasoma \cite{terasoma}, which says that
 for infinitely many points $o$ of $\mathcal{M}$ defined over a number field
$k_o$, the image of the
Galois group $Gal\,(\overline{k_o}/k_o)$ acting
on
$$H^2_{et}(A_{o,\overline{k_o}},\mathbb{Q}_l)
\oplus
H^4_{et}(X_{{o},\overline{k_o}},\mathbb{Q}_l)_{prim}$$
$$\cong
 H^2(A_{{o},\mathbb{C}},\mathbb{Q}_l)
\oplus  H^4(X_{{o},\mathbb{C}},\mathbb{Q}_l)_{prim}$$
is equal to the algebraic  monodromy group, that is, to the
image of
$\pi_{1,alg}(\mathcal{M})$ acting on
$$H^2_{et}(\mathcal{A}_{\overline{k(\mathcal{M})}},\mathbb{Q}_l)
\oplus
H^4_{et}(\mathcal{X}_{\overline{k(\mathcal{M})}},\mathbb{Q}_l).$$

As the algebraic monodromy group contains the classical monodromy group
$Im\,(\rho'_A,\rho'_U)$,
it thus follows that for such a point $o$, there is no non-zero cycle class
in $H^2_{et}(A_{o,\overline{k_o}},\mathbb{Q}_l)$
which is not proportional to
$\lambda$,
or in $H^4_{et}(X_{o,\overline{k_o}},\mathbb{Q}_l)_{prim}$, because on cycle classes,
the Galois group acts via a
finite group. As cycle classes in $H^4_{et}(X_{o,\overline{k_o}},\mathbb{Q}_l)$
 correspond bijectively via $p_*q^*$ to $(Pic\,F_{o,\overline{k_o}})
 \otimes\mathbb{Q}_l$, the proof is complete.
\cqfd
To conclude this preliminary section, recall from
\cite{voisinturin} that for any cubic $4$-fold $X$, the Fano variety
of lines $F$ carries a rational self-map
$$\phi:F\dashrightarrow F$$
which can be described as follows:

If $l\in F$ parameterizes a line $\Delta_l\subset X$, and $l$ is generic, then there
exists a unique plane $P_l\subset \mathbb{P}^5$ containing
$\Delta_l$  which is everywhere tangent to $X$ along $\Delta_l$.
As $X$ is not swept out by planes, $\Delta_l$ is not contained in any plane
contained in $X$, and thus $P_l\cap X=2\Delta_l+\Delta_{l'}$ for some point
$l'=:\phi(l)\in F$.

In the generic case, $X$ will not contain any plane, and thus the indeterminacy locus
of $\phi$
consists of those points $l\in F$ for which the plane $P_l$ is not unique,
which means equivalently that
there is a $\mathbb{P}^3\subset \mathbb{P}^5$ everywhere tangent to $X$ along $\Delta_l$.

An obvious but crucial property of $\phi$ that we will use constantly in the paper is the following
(cf \cite{voisinturin}) :
\begin{prop}\label{propertyphi}
 For  $l\in F$ away from the indeterminacy locus of $\phi$, one has the equality
$$2\Delta_l+\Delta_{\phi(l)}=h^3\,\,{\rm in}\,\,CH^3(X);$$
in particular, the class of $2\Delta_l+\Delta_{\phi(l)}$ in $CH^3(X)$ does
not depend on $l$.

\end{prop}

One first consequence is the following (see \cite{voisinturin}):
\begin{coro}  One has $\phi^*\omega=-2\omega$.
\end{coro}
Indeed, recall the incidence correspondence
$P\subset F\times X$. Then the line $\Delta_l$ has by definition its class in $CH^3(X)$ equal to
$P_*(l)$. Thus Proposition \ref{propertyphi} can be restated by saying that
\begin{eqnarray}
\label{crucialpropertyphi} P_*(\phi(l))=-2P_*l+h^3\,\,{\rm in}\,\,CH^3(X)
\end{eqnarray}
for generic $l\in F$.
Mumford's theorem conveniently generalized (cf \cite{voisinbook}, Proposition
10.24) then tells us
that we have
$$\phi^*(P^*\alpha)=-2P^*\alpha,$$
for $\alpha\in H^{3,1}(X)$ as above.
As we have $P^*\alpha=\omega$, this gives the result.
\cqfd
\section{Non preperiodicity \label{sec2}}
As mentioned in the introduction,
 Theorem  \ref{main} will be obtained as a consequence of the following
 fact :
 there exist  (many) $X$'s defined over a number field $k$, and
(many)  surfaces $\Sigma_b$ as in the previous section, also defined over a number field,
such that
the surfaces $\phi^l(\Sigma_b),\,l\in \mathbb{N},$ are Zariski dense in $X$.

Recall that we have a two-dimensional variety $B$ parameterizing very singular
 surfaces $\Sigma_b$ which are birationally equivalent to
 abelian surfaces. We shall take a desingularisation
of the total space of the family $(\Sigma_{b})_{b\in B}$. Eventually after
replacing $B$ by a Zariski open subset, this gives us a new family
parameterized  by $B$, with smooth generic fiber. The fibers
will be denoted by $\Sigma'_b$.
This section will be devoted to the proof of the following intermediate result:
Let $k(B)$ be the function field of $B$.
Consider  the following assumption on  the point $b\in B$ defined over a number field
$k(b)\subset \mathbb{C}$. We will assume $b\in B_{reg}$, the Zariski open subset
of $B$ over which the family  $(\Sigma'_{b})_{b\in B}$ is a smooth family of surfaces.
There is on one hand
the algebraic monodromy acting on the {\'e}tale cohomology of the generic  geometric fiber
$$\pi_{1,alg}(B_{reg})\rightarrow Aut\,H^1_{et}(\Sigma'_{\overline{k(B)}},\mathbb{Q}_l)\cong
Aut\,H^1_{et}(\Sigma'_{b,\overline{k(b)}},\mathbb{Q}_l) $$
and the Galois group action
acting on the {\'e}tale cohomology of the closed geometric fiber
$$Gal\,(\overline{k(b)}/k(b))\rightarrow Aut\,H^1_{et}(\Sigma'_{b,\overline{k(b)}},\mathbb{Q}_l).$$

Our assumption is :
\begin{eqnarray}\label{assumption}
{\rm  The\,\, image \,\, of\,\, the\,\, Galois\,\, group\,\, action\,\, is\,\,
equal\,\,to\,\, the\,\, algebraic\,\,
 monodromy\,\, group.}
\end{eqnarray}

By Terasoma \cite{terasoma}, this is the case for many closed points $b$ defined over a number field.

Note that a first consequence of this assumption is the fact that the
abelian surface $A_b$ is geometrically simple, and in fact has geometric
N\'eron-Severi group equal to $\mathbb{Z}$. Indeed,
we know that this is true for the abelian surface corresponding to a 
sufficiently general point of
$B(\mathbb{C})$ (see section \ref{preli}), and arguing as in the proof
of lemma \ref{conditionslemm}, this implies the same results for those fibers $A_b$
satisfying assumption (\ref{assumption}).
\begin{prop}\label{step1} Under this assumption, the
Zariski closure of the union of surfaces $\phi^l(\Sigma_b),\,l\in \mathbb{N}$,
contains at least a divisor.
\end{prop}
\begin{rema} {\rm The statement makes sense, as $\phi^l$ is always generically defined along
$\Sigma_b$. Indeed, we may assume by induction on
$l$ that $\phi^{l-1}$ is generically defined along
$\Sigma_b$. Then   $\phi^{l-1}(\Sigma_b)$ must be a surface, as proved in \cite{amerik}.
If $\phi^l$ were not
generically defined along
$\Sigma_b$, then $\phi^{l-1}(\Sigma_b)$ would be contained in the indeterminacy locus of $\phi$.
But this indeterminacy locus is a surface of general type (\cite{amerik}),
hence cannot be dominated by
a surface which is birationally equivalent to an abelian surface.}
\end{rema}
\begin{rema}{\rm  This divisor is invariant under $\phi$, hence each of its
irreducible com- ponents must be invariant under some power $\phi^l$.}
\end{rema}
Observe to start the proof of  proposition \ref{step1} that because we already know
that the $\phi^{l}(\Sigma_b)$ are surfaces, the statement
 is equivalent to the fact that there are actually
infinitely many distinct surfaces $\phi^l(\Sigma_b),\,l\in \mathbb{N}$.
Equivalently, we have to prove that we do not have
$\phi^{l'}(\Sigma_b)=\phi^{k}(\Sigma_b)$
for some $l',\,k,\,l'>k$.
This is equivalent to saying that we do not have for some $l>0,\,k\geq 0$
\begin{eqnarray}
\label{equationperiodic}\phi^{l}(\phi^{k}(\Sigma_b))=\phi^{k}(\Sigma_b),
\end{eqnarray}
or that no surface $\phi^{k}(\Sigma_b)$, $k\in \Z_{\geq 0}$, is periodic for $\phi$. In current terminology, one says that
$\Sigma_b$ is not preperiodic for $\phi$.

The proof of Proposition \ref{step1} will use the following  proposition :
\begin{prop} \label{subprop1} 1) The morphism
$$P_{b*}:CH_0(\Sigma'_{b,\mathbb{C}})\to CH_1(X_{\C})$$ induced by the pull-back to
$\Sigma'_b\rightarrow\Sigma_b\subset F$ of the incidence correspondence
$P\subset F\times X$ satisfies:
 ${P_{b*}}_{\mid CH_0(\Sigma'_{b,\mathbb{C}})_{hom}}$ factors through
the Albanese map $$CH_0(\Sigma'_{b,\mathbb{C}})_{hom}\stackrel{alb}{\rightarrow}
Alb\,\Sigma'_{b,\mathbb{C}}\cong A_{b,\mathbb{C}}.$$

2) Under assumption (\ref{assumption}), the  group morphism induced by $P_{b,*}$
$$P_{b,alb}: A_{b,\mathbb{C}}=Alb\,\Sigma'_{b,\mathbb{C}}\rightarrow CH_1(X_\mathbb{C})$$
is not of torsion.
\end{prop}
\begin{rema} {\rm It is essential here to consider the Chow groups of the
corresponding complex algebraic variety, because conjecturally
$CH_1(X)_{hom}$ is trivial for $X$ defined over a number field.}
\end{rema}
  We assume proposition \ref{subprop1}, and continue the proof of proposition
  \ref{step1}. We proceed by contradiction and assume that
  $\Sigma_{b,k}=\phi^{l}(\Sigma_{b,k})$ for some $l>1$, where
  $\Sigma_{b,k}:=\phi^k(\Sigma_b)$.
Let us choose a desingularization
$\Sigma'_{b,k}$ of $\Sigma_{b,k}$. We first prove:
  \begin{lemm}\label{alksigmabk} The Albanese map
  $$alb: \Sigma'_{b,k}\rightarrow Alb\,{\Sigma'_{b,k}}$$
  is a birational map.
\end{lemm}
{\bf Proof.} First observe that $\Sigma'_{b,k}$ is rationally dominated via the rational
map $\phi^k$ by
the abelian surface $A_b=Alb\,{\Sigma'_b}$. On the other hand, we claim that
$Alb\,{\Sigma'_{b,k}}$ is isogenous (via $\phi^k$) to $Alb\,{\Sigma'_{b}}$.
Indeed, this follows from the following three facts:
We know from proposition \ref{subprop1}, 1)  that for $\Sigma'_b$ (hence
 also for $\Sigma'_{b,k}$ by property (\ref{crucialpropertyphi})), the
 restriction $P_{b*}:CH_0(\Sigma'_{b,\mathbb{C}})_{hom}\rightarrow CH_1(X_\mathbb{C})$ (resp.
 $P_{b,k,*}:CH_0(\Sigma'_{b,k,\mathbb{C}})_{hom}\rightarrow CH_1(X_\mathbb{C})$) factors
 through $Alb\,\Sigma'_{b,\mathbb{C}}$ (resp. $Alb\,\Sigma'_{b,k,\mathbb{C}}$).  Let
 $$P_{b,k,alb}:Alb\,{\Sigma'_{b,k,\mathbb{C}}}\rightarrow CH_1(X_\mathbb{C}),\,\,
 P_{b,alb}: A_{b,\mathbb{C}}\rightarrow CH_1(X_\mathbb{C})$$
be induced respectively by the maps $P_{b,k,*}$ and $P_{b*}$,
  and let
  $$\phi^k_*:A_{b,\mathbb{C}}=Alb\,\Sigma'_{b,\mathbb{C}}\rightarrow Alb\,{\Sigma'_{b,k,\mathbb{C}}}$$
  be  induced by the rational map $\phi^k$ between $\Sigma'_b$ and
  $\Sigma'_{b,k}$.
 By formula (\ref{crucialpropertyphi}), factoring everything through the Albanese
 maps,  we get
\begin{eqnarray}\label{memeformule}P_{b,k,alb}\circ \phi^k_*=(-2)^kP_{b,alb}:A_{b,\mathbb{C}}
\rightarrow CH_1(X_\mathbb{C}).
\end{eqnarray}

On the other hand,
we know by Proposition \ref{subprop1}, 2) that $P_{b, alb}: A_{b,\mathbb{C}}
\rightarrow CH_1(X_\mathbb{C})$ is not of torsion.
As $P_{b,*}$ is a group morphism which is induced by a correspondence, its kernel is
a countable union of translates of an abelian subvariety of $A_{b,\mathbb{C}}$. Thus it must
be countable
because we know  that $A_b$ is a geometrically simple abelian variety
and that $P_{b,*}$ is non zero.
Formula (\ref{memeformule}) then shows that the surjective morphism
$\phi^k_*$ has  countable kernel, hence that it is an isogeny.

In conclusion, the surface $\Sigma'_{b,k}$ sits between two abelian surfaces:
$$ A_b\dashrightarrow \Sigma'_{b,k}\stackrel{alb}{\rightarrow} Alb\,\Sigma'_{b,k},$$
where the composed map is an isogeny. It follows immediately that
$alb: \Sigma'_{b,k}{\rightarrow} Alb\,{\Sigma'_{b,k}}$ is birational.
\cqfd
As a consequence of this lemma and formula (\ref{memeformule}),
we get now:
\begin{coro}\label{corointer} The rational map
$$\phi^l:\Sigma'_{b,k}\rightarrow \Sigma'_{b,k}$$
has degree $(-2)^{4l}=(16)^l$.
\end{coro}
{\bf Proof.}  Indeed,
formula (\ref{memeformule}) applied to the maps $\phi^k$ and $\phi^{k+l}$ tells us that
\begin{eqnarray}
\label{above15juillet}
P_{b,k,alb}\circ\phi^l_*=(-2)^lP_{b,k,alb}: Alb\,{\Sigma'_{b,k,\mathbb{C}}}\rightarrow CH_1(X_\mathbb{C}).
\end{eqnarray}
Here $\phi^l_*:Alb\,{\Sigma'_{b,k,\mathbb{C}}}\rightarrow Alb\,{\Sigma'_{b,k,\mathbb{C}}}$
 is induced by the self-rational map
$\phi^l:{\Sigma'_{b,k}}\rightarrow {\Sigma'_{b,k}}$.
Furthermore,  we know that
$Ker\,P_{b,k,alb}$ is countable.
 From  formula (\ref{above15juillet}), we
 thus  deduce that $\phi^l_*$ must be multiplication by $(-2)^l$, hence must have degree
 $(-2)^{4l}$.
As  the  surface
 $\Sigma'_{b,k}$ is birationally equivalent to its Albanese variety by lemma
 \ref{alksigmabk}, the same result holds for
 $\phi^l:\Sigma'_{b,k}\rightarrow \Sigma'_{b,k}$.

\cqfd
{\bf Proof of Proposition \ref{step1}}.
We know that $\phi^l$ is generically well defined along $\Sigma'_{b,k}$.
Take the standard desingularization $\widetilde{\phi}:\widetilde{F}\rightarrow F$
with $\pi: \widetilde{F}\rightarrow F$ the blow-up of indeterminacy locus
of $\phi$. It is proved in \cite{amerik} that
$\widetilde{\phi}$ does not contract divisors and can only contract a
surface which is already contracted by $\pi$ (in fact even this is impossible,
but needs some extra argument). It follows that for
the natural  desingularization
$\widetilde{\phi^l}:\widetilde{\widetilde{F}}\rightarrow F$,
$\pi':\widetilde{\widetilde{F}}\rightarrow F$ of $\phi^l$
obtained by taking
the graph, the following holds:
For any  subvariety $\Sigma\subset\widetilde{\widetilde{F}} $
which is contracted under
$\widetilde{\phi^l}$,  the image of
$\Sigma$ under $\pi'$ is at most $1$-dimensional.
We apply this to any component $\Sigma$ of
$\widetilde{\phi^l}^{-1}(\Sigma_{b,k})$.  We conclude that
if either $dim\,\Sigma\geq 3$ or $\Sigma$ does not dominate $\Sigma_{b,k}$ via
$\widetilde{\phi^l}$, $\pi'(\Sigma)$ has dimension at most $1$.
On the other hand,
$\pi'(\widetilde{\phi^l}^{-1}(\Sigma_{b,k}))$ contains $\Sigma_{b,k}$ by periodicity, and we
have just seen that the restriction of $\phi^l$ to $\Sigma_{b,k}$  is
already of degree $(16)^l=deg\,\widetilde{\phi^l}$. It thus follows that   the only
$2$-dimensional component of
$\widetilde{\phi^l}^{-1}(\Sigma_{b,k})$ dominating $\Sigma_{b,k}$ is
sent birationally by $\pi'$ to
$\Sigma_{b,k}$, while the non dominating components or the three dimensional
component are contracted to curves or points in
$F$ by $\pi'$.

 This implies  that we have the equality
of $2$-dimensional cycles:
$$(\phi^l)^{*}(\Sigma_{b,k})=\Sigma_{b,k}$$ and thus
$$(\phi^l)^{*}([\Sigma_{b,k}])=[\Sigma_{b,k}],$$
where $[\,\cdot\,]$ denotes the cohomology class of a cycle.
This contradicts the computation of the eigenvalues of $(\phi^l)^*$ on $H^4(F_\mathbb{C},\mathbb{Q})$
performed in \cite{amerik}.

The proof is now complete, assuming Proposition \ref{subprop1}.
\cqfd
{\bf Proof of Proposition \ref{subprop1}.}
Statement 1) follows from the fact that the surface $\Sigma_b$ is the surface of lines of a
nodal cubic
threefold $Y_b$, for which we know that the Abel-Jacobi map on a desingularization
$\widetilde{Y}_b$ induces an isomorphism:
$$CH_1(\widetilde{Y}_{b,\mathbb{C}})_{hom}\stackrel{\Phi_{\widetilde{Y}_b}}{\cong}
J^3(\widetilde{Y}_{b,\mathbb{C}}).$$
The correspondence from $\Sigma'_b$ to $X$ thus factors to a correspondence from
$\Sigma'_b$ to $\widetilde{Y}_b$, which, when restricted to
$CH_0(\Sigma'_{b,\mathbb{C}})_{hom}$ has to factor through $Alb\,\Sigma'_{b,\mathbb{C}}=J^3(\widetilde{Y}_{b,\mathbb{C}})$.

For the proof of statement 2), we will follow the argument of
\cite{ggp}. In fact this will be   as in \cite{ggp}
a consequence  of  an even stronger
statement (Proposition \ref{plusfort}), the proof of which is attributed
there to Bloch and Nori, and  which we will also need
in section \ref{conclusion}.

Let us first introduce a corrected codimension $3$ cycle
$$P'\in CH^3( F\times X)_\mathbb{Q}$$
which has the property that its cohomology class satisfies
\begin{eqnarray}
\label{formPprime}[P']\in H^2(F_\mathbb{C},\mathbb{Q})\otimes H^4(X_\mathbb{C},\mathbb{Q})_{prim},
\end{eqnarray}
and which differs from the incidence cycle $P$ by a combination of cycles of the form
$\Gamma_i\times h^i$, $\Gamma_i\in CH^{3-i}(F)$.
This cycle is obtained by considering
the  cycle $j_*P$ on $F\times \mathbb{P}^5$, where $j$ is the inclusion of $X$
in $\mathbb{P}^5$ (or of $F\times X$ in $X\times \mathbb{P}^5$). Denoting by
$H\in CH^1(\mathbb{P}^5)$ the class of $\mathcal{O}_{\mathbb{P}^5}(1)$, this cycle
$j_*P\in CH^4(F\times
\mathbb{P}^5)$ can be written as $\sum_{i=1}^{i=4}pr_1^*\Gamma_i\cdot pr_2^*H^{i}$, with $\Gamma_i\in
CH^{4-i}(F)$ (cf \cite{voisinbook}, Theorem 9.25).
The fact that the sum is only over $i\geq1$ is due to the fact that
$H^5 j_*P=0$.
Let $$P':=P-\frac{1}{3}(\sum_{i=0}^{i=3}pr_1^*\Gamma_{i+1}\cdot pr_2^*h^{i})\in CH^3(F\times X)_\mathbb{Q},$$
where now the projectors are those of $F\times X$ to its factors.
Then $j_*P'=0$ because $j_*h^i=3H^{i+1}$, and property
(\ref{formPprime})  follows by Lefschetz theorem on hyperplane sections, which tells that
$Ker\,(j_*:H^*(X_\mathbb{C},\mathbb{Q})\rightarrow H^{*+2}(\mathbb{P}^5_\mathbb{C},\mathbb{Q}))$ is zero for
$*\not=4$.

We can see, using Poincar{\'e} duality on $X$, the vector space
$H^2(F_\mathbb{C},\mathbb{Q})\otimes H^4(X_\mathbb{C},\mathbb{Q})_{prim}$ as
$$Hom\,(H^4(X_\mathbb{C},\mathbb{Q})_{prim},H^2(F_\mathbb{C},\mathbb{Q})).$$ Thus,
for any morphism $f: T\rightarrow F$, the cycle
$f^*P'$ is cohomologous to $0$ if and only if the map
$$f^*\circ P^*:H^4(X_\mathbb{C},\mathbb{Q})_{prim}\rightarrow H^2(T_\mathbb{C},\mathbb{Q})$$
is $0$.  (Note that, by definition of $P'$, $P^*={P'}^*$ on $H^4(X_\mathbb{C},\mathbb{Q})_{prim}$.)

Assume now  that this last condition is satisfied. (This is true in particular if
$T$ is one of the surfaces $\Sigma'_b$, by the same argument
as in section \ref{preli}, where we prove that the surfaces $\Sigma_b$ are Lagrangian with respect to
$\omega$.)

Then the cycle $f^*P'$ is cohomologous to $0$. If furthermore
$T$ and $f$ are defined over a number field $k$, then by the comparison theorems
between {\'e}tale $l$-adic and Betti cohomology for algebraically closed fields of characteristic $0$,
the cycle $f^*P'$ has trivial  $l$-adic {\'e}tale cohomology class in
$H^6_{et}(T_{\overline{k}}\times_{\overline{k}} X_{\overline{k}},\mathbb{Q}_l(3))$.
On the other hand, the cycle
$f^*P'$ defined over  $k$ has a {\it continuous {\'e}tale} (cf \cite{jannsen}, \cite{raskind})
cohomology class
$$[f^*P']\in
H^6(T_{k}\times_{{k}} X_{{k}},\mathbb{Q}_l(3)).$$

Following \cite{raskind}, we omit in the sequel both ``continuous'' and ``\'etale'' in the notation
for these cycle classes.

The Hochschild-Serre spectral sequence for continuous {\'e}tale cohomology then gives us
the $l$-adic Abel-Jacobi invariant
$$\gamma_{f^*P'}\in H^1(k, H^5
(T_{\overline{k}}\times_{\overline{k}} X_{\overline{k}},\mathbb{Q}_l(3))),$$
which lies in fact (because the cycle $P'$ is deduced from $P$ by
a Chow-K\"unneth decomposition) in
the subspace
$H^1(k, H^1(T_{\overline{k}},\mathbb{Q}_l)
\otimes_{\mathbb{Q}_l} H^4(X_{\overline{k}},\mathbb{Q}_l)(3))$.

(Here $H^1(k,\cdot)$ denotes the continuous \'etale cohomology of $Spec\,k$.)

Coming back to the case where $T$ is one of our surfaces $\Sigma'_b$, with $b\in B$ a point
defined over a number field $k(b)$, let us denote by
$P'_b\in CH^3(\Sigma'_b\times X)_\mathbb{Q}$ the cycle $\tau_b^*P'$, where
$\tau_b:\Sigma'_b\times X\rightarrow F\times X$
is the natural map.
Using
 the Bloch-Srinivas lemma (cf \cite{blS}
or \cite{voisinbook}, Corollary 10.20), we get as in \cite{ggp} the following
lemma:
\begin{lemm} \label{intermediate} If the $l$-adic Abel-Jacobi invariant
$$\gamma_{P'_b}\in H^1(k(b), H^1(\Sigma'_{\overline{k(b)}},\mathbb{Q}_l)
\otimes_{\mathbb{Q}_l} H^4(X_{\overline{k(b)}},\mathbb{Q}_l)(3))$$
is not of torsion, then the morphism
$$P_{b,alb}: Alb\,\Sigma'_{b,\mathbb{C}}=A_{b,\mathbb{C}}\rightarrow CH_1(X_\mathbb{C})$$
is not of torsion.
\end{lemm}

Using this lemma, the statement 2) of proposition \ref{subprop1} is a consequence of the following
proposition:

\begin{prop} \label{plusfort} Under the assumption (\ref{assumption}), the
$l$-adic Abel-Jacobi invariant
$$\gamma_{P'_b}\in H^1(k(b)), H^1(\Sigma'_{\overline{k(b)}},\mathbb{Q}_l)
\otimes_{\mathbb{Q}_l} H^4(X_{\overline{k(b)}},\mathbb{Q}_l)(3))$$
is non-zero.
\end{prop}
{\bf Proof.}
The proof of this statement goes as in \cite{ggp}, once we notice the following facts:
We have a $2$-dimensional family of surfaces which sweep out $F$ :
$$\mathcal{S}_B:=\cup_{b\in B}\Sigma_b'\stackrel{s}{\rightarrow} F$$
 where $s$ is dominating. For a generic curve $B'\subset B$
defined over a number field, the family of surfaces
$$\mathcal{S}_{B'}:=\cup_{b\in B'}\Sigma_b'\stackrel{s}{\rightarrow} F$$
covers an hypersurface in $F$, and it thus follows that
the pull-back map by $s$:
$$ H^{2,0}(F)\rightarrow H^{2,0}(\mathcal{S}_{B'})$$
is non-zero because the $(2,0)$-form $\omega$ is non-degenerate, hence cannot vanish
on a divisor of $F$.

More precisely, for any Zariski open set $B'_0\subset B'$, the cohomology class
$[s^*\omega]$ does not vanish in the Betti cohomology space
$H^2(\mathcal{S}_{B'_0,\mathbb{C}}, \mathbb{C})$, because it is of type $(2,0)$.
We will take for $B'_0$ the locus where the fiber of our family of surfaces
$\pi:\mathcal{S}_{B'}\rightarrow B'$ is smooth.

Hence we conclude similarly that the Betti
cohomology class of the cycle
$$s^*P'_{\mid \mathcal{S}_{B'_0}}\in CH^3(\mathcal{S}_{B'_0}\times X)$$
does not vanish in
$$ H^6(\mathcal{S}_{B'_0,\mathbb{C}}\times X_\mathbb{C},\mathbb{Q}).$$
However, we know already that the cycle $s^*P'$ restricted to the fibers
$\Sigma'_b\times X$ of $\pi$ has vanishing cohomology class. Hence
$[s^*P'_{\mid \mathcal{S}_{B'_0}}]$ lies at least in the $L^1$ level of the Leray filtration
associated to
$$\pi'=\pi\circ pr_1:\mathcal{S}_{B'_0,\mathbb{C}}\times X_\mathbb{C}\rightarrow B'_{0,\mathbb{C}},$$
 and gives an element of
$$H^1(B'_{0,\mathbb{C}},R^5\pi'_{\mathbb{C}*}\mathbb{Q})$$
and in fact more precisely (using the fact that $P'$ was obtained
from $P$ by applying a Chow-K\"unneth projector on
$X$) in
$$H^1(B'_{0,\mathbb{C}},R^1\pi_{\mathbb{C}*}\mathbb{Q}\otimes H^4(X_\mathbb{C},\mathbb{Q})_{prim}).$$
This last class is non-zero because the Leray spectral sequence has only two terms, as
the affine curve $B'_0$ has the homotopy type of a CW-complex of dimension $1$.

We can now make the same construction in $l$-adic {\'e}tale cohomology for
the corresponding algebraic varieties defined over $\overline{\mathbb{Q}}$.
Using the comparision theorems between Betti cohomology of the complex manifold and {\'e}tale cohomology
of the algebraic variety, we conclude that
the {\'e}tale cohomology class
$$[p^*P'_{\mid \mathcal{S}_{B'_0}}]_{et}\in H^6_{et}
(\mathcal{S}_{B'_0,\overline{\mathbb{Q}}}\times_{\overline{\mathbb{Q}}} X_{\overline{\mathbb{Q}}},
\mathbb{Q}_l(3))$$
is non zero, lies in the $L^1$ level of the Leray filtration for {\'e}tale
cohomology, and has a non zero image in
$$H^1_{et}(B'_{0,\overline{\mathbb{Q}}},R^1\pi_{et*}\mathbb{Q}_l\otimes_{\mathbb{Q}_l}
 H^4(X_{\overline{\mathbb{Q}}},\mathbb{Q}_l)_{prim}(3)).$$

On the other hand, let us  choose a point $b\in B'_0(\overline{\mathbb{Q}})$ satisfying the
assumption (2.2).
Then we conclude as in \cite{ggp}, Lemma 7, from the non vanishing
above,  that the class
 $$\gamma_{P'_b}\in H^1(k(b), H^1(\Sigma'_{\overline{k(b)}},\mathbb{Q}_l)
\otimes_{\mathbb{Q}_l} H^4(X_{\overline{k(b)}},\mathbb{Q}_l)(3))$$
is non-zero.

\cqfd
\section{Study of invariant divisors}
Our goal in this section is to exclude  the case where
the Zariski closure of the union of the surfaces
$\Sigma_{b,i}$ is a divisor.

Reasoning by contradiction,
we will denote in the sequel by $\tau:D\rightarrow F$ a desingularization of a
$3$-dimensional  irreducible
component of the Zariski closure of the union of the surfaces $\Sigma_{b,i}$, where $b$ is fixed and is supposed to
satisfy the conclusion of Proposition \ref{step1}.
As this Zariski closure has only finitely many $3$-dimensional components and is stable under
$\phi$, it follows that each of its irreducible components is stable under some power
$\phi^l$ of $\phi$.

We will denote by
 $\Sigma_{b,n_i}$ the surfaces $\Sigma_{b,k}$ contained
in $\tau(D)$ and  $\Sigma'_{b,n_i}$ their proper transforms under $\tau$.

We shall also denote by $\omega_D$ the holomorphic $2$-form
$\tau^*\omega$. Notice that $\omega_D$ is a non-zero form, because $\omega$ is non degenerate.
Thus at a point where $\tau $ is of maximal rank $3$, $\omega_D$ cannot vanish.

\subsection{A Jouanolou type argument \label{secclassi}}

The kernel of $\omega_D$ at the generic point is one-dimensional. This
provides a saturated
invertible subsheaf $L\subset T_D$, generically defined as
$Ker\,\omega_D:T_D\rightarrow\Omega_D$. As the surfaces $\Sigma'_{b,n_i}$ are
isotropic for $\omega_D$, $L$ must be tangent to any
$\Sigma'_{b,i}$, or at least to those not contained in the singular set of
$\tau(D)$. Recall that the  surfaces $\Sigma'_{b,n_i}$ are Zariski dense
in $D$. In particular there are infinitely many such surfaces which are
distinct and not contained in the singular set of
$\tau(D)$.

 We start by proving an analogue of
Jouanolou's theorem (cf \cite{J}) for this situation:

\begin{prop}\label{jouan} There exists a rational map
$f: D\dasharrow \p^1$ such that $L$ is tangent to the fibers of $f$.
\end{prop}
{\bf Proof:} For convenience of the reader, we recall at the same time
the construction
of Jouanolou. For large $N$, the surfaces $\Sigma'_{b,n_i},\,i\leq N$ are not linearly
independent in $NS(D)$. Therefore one can find integers $m_i$ such that
the line bundle associated to the divisor $M=\sum_{i=i_0}^Nm_i\Sigma'_{b,n_i}$
has zero Chern class in $H^1(D, \Omega_D)$. Here, $i_0$ is chosen in such way that
 $\omega_D$ does not vanish on $\Sigma'_{b,n_i}$ for $i\geq i_0$.
Let $U_{\alpha}$ be an open
covering of $D$, and $f_{\alpha}$ be a meromorphic function defining $M$ on $U_{\alpha}$. The vanishing
of the first Chern class of the cocycle $g_{\alpha \beta}\in
H^1(D, {\cal O}^*_D)$ means that on the intersections
$U_{\alpha\beta}=U_{\alpha}\cap U_{\beta}$
$$\frac{df_{\alpha}}{f_{\alpha}}-\frac{df_{\beta}}{f_{\beta}}=
\sigma_{\alpha}-\sigma_{\beta}$$
for some $\sigma_{\alpha}\in \Omega^1_{U_{\alpha}}$,
$\sigma_{\beta}\in \Omega^1_{U_{\beta}}$. Therefore, the
$\frac{df_{\alpha}}{f_{\alpha}}-\sigma_{\alpha}$ patch together into a
meromorphic 1-form $\psi_M$ with logarithmic singularities
along the components of the support of
$M$, defined up to a holomorphic 1-form, that is, we have
a map
$$\psi: Div^0_{\Sigma}(D)\otimes \C \to H^0(D, \Omega^1_D
\otimes {\cal M}^*_D)/H^0(D,\Omega^1_D),$$
where ${\cal M}^*_D$ denotes the sheaf of meromorphic functions and
where $Div^0_{\Sigma}(D)$ is the kernel of
the cycle class map, restricted to the subspace generated by the $\Sigma'_{b,n_i},\,i\geq i_0$.
 Moreover,
on $U_{\alpha \beta}$, $d\sigma_{\alpha}=d\sigma_{\beta}$, so there is
a well-defined holomorphic 2-form $d\sigma$. It follows from the
degeneration of ``Hodge to de Rham" spectral sequence that $d\sigma$
is actually zero (cf. \cite{J}). Therefore $\psi_M$ is closed.

Moreover, the map $\psi$ is injective (this is \cite{J}, Lemma 2.8).

The meromorphic form $\psi_M$ induces a meromorphic section of $L^*$.
Since outside of the singularities of the foliation, that is, outside of an
analytic subset of codimension two in $D$, $L$ is tangent to $\Sigma'_{b,n_i}$
along $\Sigma'_{b,n_i}$ for $i\geq i_0$,
this section is actually holomorphic, because $\psi_M$ has logarithmic
singularities along the smooth part of the  $\Sigma'_{b,n_i}$'s appearing in the support of $M$.
So the map $\psi$ induces a map
$$\chi:  Div^0_{\Sigma}(D)\otimes \C \to H^0(D,L^*)/V,$$
where $V$ is the image in $H^0(D,L^*)$ of the space of global holomorphic 1-forms on $D$.
The kernel of $\chi$ is clearly infinite-dimensional.

Let $M$ and $M'$
be two non-proportional elements
of  $Ker\,\chi\subset Div^0_{\Sigma}(D)\otimes \C$.
Then we can choose the meromorphic forms $\psi_M$ and $\psi_{M'}$ in such a way that
$$\psi_{M\mid L}=0,\,\,\psi_{M'\mid L}=0.$$
It then follows that $L$ is contained in the kernel
of the meromorphic $2$-form $\psi_M\wedge \psi_{M'}$. But this property is also satisfied
by our  original 2-form
$\omega_D$. Therefore $\psi_M\wedge \psi_{M'}=f\omega_D$ for some meromorphic function
$f$ on $D$.
We meet now two cases:

1) Assume that  $f$ is always a constant function on $D$.
In this case, the meromorphic form $\psi_M\wedge \psi_{M'}$ is nonsingular.
On the other hand, the meromorphic $1$-form $\psi_{M'}$ has logarithmic singularities
along the surfaces $\Sigma'_{b,n_j}$ appearing in $M'$. The absence of singularities
in $\psi_M\wedge \psi_{M'}$ then implies that the restriction of $\psi_M$
to the surfaces $\Sigma'_{b,n_j}$ appearing in $M'$ and not in $M$ is zero.
But as $M'$ is arbitrary, it then follows that there are infinitely
many surfaces $\Sigma'_{b,n_i}$ which are integral surfaces for the $1$-form
$\psi_M$. Jouanolou's theorem
\cite{J} then tells that there exists a non constant
meromorphic  function from
 $D$ to $\mathbb{P}^1$,  whose fibers are leaves of $\psi_M$. As $\psi_M$ vanishes on $L$,
 we are done in this case.

2) In the second case,  we have for some
$M'$ a non-constant rational map $f: D\dasharrow \p^1$.
Finally, differentiating
the equality $\psi_M\wedge \psi_{M'}=f\omega_D$,
an using the fact  that all our forms are closed, we get $df\wedge\omega_D=0$, so $df$ vanishes on
$Ker\,\omega_D=L$.

\cqfd

\begin{coro}\label{twofibrat}
We have two possibilities: either

(1) $D$ is rationally fibered in (birationally) abelian varieties,
and $\Sigma'_{b,n_i}$
are fibers for all but finitely many $i$'s; or

(2) $D$ is rationally fibered in curves over a surface $T$, the fibers are
integral curves of the foliation $L$, and $\Sigma'_{b,n_i}$
project onto curves in $T$. Furthermore the $2$-form $\omega_D$ is pulled-back from a
$2$-form on $T$.
\end{coro}
Here,``rationally fibered" means that the fibration map is only rational,
not necessarily regular. We can, of course, arrange for it to be regular
by blowing $D$ up.

\medskip

\noindent{\bf Proof:} Consider the Stein factorisation
$g: D\dasharrow C$ of the map $f$. If all but finitely many
$\Sigma'_{b,n_i}$'s are fibers, we are done; if not, consider the general
fiber $D_c$. On $D_c$, we still have the foliation $L$ from the
proposition \ref{jouan}. This foliation has infinitely many integral curves, which
are components of $\Sigma'_{b,n_i}\cap D_c$, where $i$ is such that
$\Sigma'_{b,n_i}$ is not a fiber. So, by the original
Jouanolou's theorem, there must be a fibration
$h_c:D_c\dasharrow Z_c$ over a curve, whose fibers are generically connected
and tangent to
$L$. A countability argument for the Chow varieties
of curves in $D$  shows that we can assume
the   fibers of $h_c$, $c\in C$, form a family of algebraic curves covering $D$ and consisting of
integral leaves of $L$.
It is then immediate that these curves provide a fibration of $D$ to a surface
$T$, which we may assume to have connected fibers.

As $L$ is tangent to these fibers, and to the surfaces $\Sigma'_{b,n_i}$, it follows that the
$\Sigma'_{b,n_i}$'s
 are contracted to curves in $T$. Finally, as the fibration has connected fibers, whose tangent space
 is generically the kernel of $\omega_D$, it follows that $\omega_D$ comes from a
 form on $T$.

\cqfd

\subsection{Case where $D$ is fibered into abelian surfaces \label{secclassigeo}}

In this subsection, we will exclude case 1) of Corollary \ref{twofibrat}. We argue again by contradiction,
and assume that we have a divisor $D$, invariant under some power $\phi^l$ 
of $\phi$, and
admitting a fibration over a curve, whose fibers are surfaces birationally equivalent to abelian surfaces.
We know furthermore that infinitely many fibers are surfaces $\Sigma'_{b,n_i}$, which
are not periodic under $\phi^l$.

The starting point is the following

\begin{prop}\label{kod}
The Kodaira dimension of $D$ is zero.
\end{prop}
{\bf Proof:} $D$ is (rationally) fibered over a curve in varieties
of Kodaira dimension zero. By the universal property of the Iitaka
fibration, $\kappa(D)$ is thus at most one.

To exclude the case where $\kappa(D)=1$, we use  the result of \cite{NZ}, saying that
a rational self-map $\psi$ of any variety $X$ induces an automorphism
of finite order on the base $B$ of the Iitaka fibration (note that the fact
that $\psi$ induces a polarization-preserving automorphism of $B$ is clear
because the pull-back $\psi^*$
is an automorphism of the vector space of pluricanonical forms; the non-trivial
statement is the finiteness of the order of the induced automorphism  on $B$. Also, we do not need the full strength of
Theorem A of \cite{NZ}; the crucial Proposition 2.3 for an abelian fibration
suffices, and this is shown by a monodromy argument).

 We now apply this to $X=D$ and $\psi=\phi^l$. If $\kappa(D)=1$,
then the map $f:D\dasharrow \p^1$ constructed in \ref{jouan} is the Iitaka
fibration, and the finiteness of the order of the induced map on $\p^1$ implies
that the $\Sigma_{b,n_i}$ (which are fibers of $f$) are periodic, whereas we know that this is not the case.

Suppose now that $\kappa(D)=-\infty$. Then $D$ is uniruled. Consider the
rationally connected  fibration $r:D\dasharrow Q$. As the fibers are rationally connected, every
holomorphic form on $D$ must be pulled-back from $Q$. Since $D$ carries a non-zero
holomorphic 2-form $\omega_D$, $Q$ is a surface, and $\omega_D=r^*\eta$
for some $0\not=\eta\in H^{2,0}(Q)$.
The surfaces $\Sigma'_{b,n_i}$ cannot  project onto curves on $Q$, since
they are not uniruled. Therefore they dominate $Q$ via $r$. But as they are
isotropic for $\omega_D$, we get
$$\omega_{D\mid \Sigma'_{b,n_i}}=0=r^*\eta_{\mid \Sigma'_{b,n_i}},$$ hence
$\eta=0$, which is  a contradiction.

\cqfd

If $D$ had numerically trivial canonical class, the natural idea to continue
would be to use the Bogomolov-Beauville decomposition. In our situation,
we can take the minimal model of $D$, but it can have terminal singularities,
and at this moment, no analogue of Bogomolov-Beauville decomposition is known
for terminal varieties. Nevertheless, it exists for threefolds with
a holomorphic $2$-form.

The following proposition is proved in \cite{CP}; we give a somewhat more
detailed argument for reader's convenience.

\begin{prop}\label{mmp} (Theorem 3.2 of \cite{CP})
For any smooth projective threefold $D$ such that $\kappa(D)=0$ and
$h^{2,0}(D)\neq 0$,
there exists a dominant rational map $\mu:Y\dasharrow D$, where $Y$ is either
an abelian threefold, or a product $E\times S$ of a K3 surface $S$
and an elliptic curve $E$.

\end{prop}

\medskip

\noindent{\bf Proof:} By minimal model theory, there is a birational map
$D\dasharrow D^0$, such that $D^0$ has at most canonical singularities,
the Weil divisor $K_{D^0}$ is $\Q$-Cartier and $mK_{D^0}=0$ for
some $m>0$. Take the canonical cover $D'=Spec(\oplus_{i=0}^{m-1}{\cal O}
(iK_{D^0}))$ (where ${\cal O}(iK_{D^0})$ are sheaves associated to Weil
divisors; see \cite{R}, appendix to section 1, for a discussion of the
definition of the pluricanonical sheaves etc. in this case). It is well-known
(cf \cite{R}) that it is \'etale in codimension 2 (wherever
$K_{D^0}$ is Cartier), the singularities of $D'$
 are canonical and that $K_{D'}=0$. Let $Y$ be a resolution of $D'$.

Since canonical singularities are rational by \cite{E},
 that is, there is no higher direct images of the structure sheaf,
we have $H^{2,0}(Y)\neq 0$ (and $H^0(D', \Omega^2_{D'})\neq 0$ for a
reasonable  definition of $\Omega^2_{D'}$). The argument of Peternell
(\cite{P}, section 5), generalization of that of Bogomolov's to the
singular case,
 produces, for a form $\sigma \in H^0(Y,\Omega^2_Y)$,
a holomorphic 1-form $\eta \in  H^0(Y,\Omega^1_Y)$ such that
$\eta\wedge\sigma$ generates $H^0(Y, K_Y)$.

The Albanese map of $Y$ factors through $D'$ and gives a fibration
$\alpha:D'\to Alb\,D':= Alb\,Y$ (\cite{Kaw}). Now the results of
\cite{Kaw}, section 8, say that this becomes a product after a finite \'etale
covering. Repeating eventually the process (i.e. considering eventually
the Albanese map) for the fiber of $\alpha$, we get the result.

\cqfd

Let us denote by $H^2_0(F, \Q)$ the orthogonal of
$c_1(\mathcal{L})$ with respect to the Beauville-Bogomolov form. This is just
the  image of the map
$p_*q^*:H^4(X,\mathbb{Q})_{prim}\rightarrow H^2(F, \Q)$ (cf section \ref{preli}).

\begin{lemm}\label{hodge}
For $F$ such that $Pic(F)=\Z$, $H^2_0(F, \Q)=H^2_{tr}(F, \Q)$ is a
simple
Hodge
structure.
\end{lemm}
{\bf Proof:} This follows from $h^{2,0}(F)=1$ by a standard argument:
let $V$ be a non trivial Hodge
substructure of $H^2_0(F, \Q)$. If $V$ is of level zero, that is,
$V\otimes \C\subset H^{1,1}(F)$, then there are integral $(1,1)$ classes
in $H^2_0(F, \C)$, contradicting $Pic(F)=\Z$. If $V$ is not of level zero,
$V\otimes \C$ contains $H^{2,0}(F)$, so the orthogonal of $V$ is of level
zero, hence must be trivial, and we conclude that $V=H^2_0(F, \Q)$.

\cqfd

For a projective manifold $X$, denote by $b_2^{tr}(X)$ the difference
between $b_2(X)$ and the Picard number $\rho(X)$.

\begin{lemm}\label{betti}
Let $Y,\,Z$ be smooth projective of the same dimension. If there exists
a dominant rational map $f:Y\dasharrow Z$, then $b_2^{tr}(Z)\leq b_2^{tr}(Y)$.
\end{lemm}
{\bf Proof:} Let $p:Y'\to Y$ be a resolution of $f$, so that $p$ is a
composition
of blow-ups with smooth centers and $q:Y'\to Z$ is a generically finite
morphism, $q=f\circ p$. Since $q_*q^*=deg(q)\cdot Id$ on  cohomology groups,
$q^*:H^2(Z,\C)\to H^2(Y',\C)$
is injective and takes transcendental cohomology to transcendental cohomology. So
$b_2^{tr}(Z)\leq b_2^{tr}(Y')$. Since $p$ is a composition of blow-ups,
$b_2^{tr}(Y')=b_2^{tr}(Y)$.

\cqfd

We now conclude as follows. From propositions
\ref{kod}  and \ref{mmp}, our divisor
$D$ is rationally dominated either by an abelian threefold, or by a
product of an elliptic curve and a K3-surface. Consider
the morphism of Hodge structures
$$\tau^*: H^2_{tr}(F, \C)\to H^2(D, \C)$$
As  $\omega$ restricts
non-trivially to $D$,  this morphism is non-zero. By irreducibility
of the Hodge structure on $H^2_{tr}(F, \Q)$, it is injective and its image
is contained in  the transcendental part of $H^2(D,\mathbb{Q})$,
so $b_2^{tr}(D)$ is at least
$22$ (recall that $b_2(F)=23$ since $F$ is deformation equivalent to
the $Hilb^2$ of a K3 surface \cite{be-do} and $F$ is generic, so
the transcendental dimension is only one less).
This is clearly greater than $b_2^{tr}$ of a three-dimensional torus.
Finally, let $E$ be an elliptic curve and $S$ a K3-surface; we have
$$H^2(E\times S,\Q)=H^2(E,\Q)\otimes H^0(S,\Q)\oplus H^0(E,\Q)
\otimes H^2(S,\Q);$$
this is of dimension $23$ and has at least a $2$-dimensional subspace of
algebraic classes, so that $b_2^{tr}(E\times S)$ is again strictly less
than $22$, contradicting lemma \ref{betti}.

We thus proved  that the case of a fibration in abelian surfaces is impossible.

\subsection{Case where $D$ is fibered over a surface \label{conclusion}}
In this section we consider the remaining situation, where $D$ admits a rational fibration
 over a surface $T$, in such a way that the $2$-form $\omega_D$ is pulled-back from
 a $2$-form on $T$
 and the
 countably many Zariski dense surfaces $\Sigma'_{b,n_i}\subset D$ contained in $D$
  are inverse images of curves
 $T_i\subset T$:
 $$\Sigma'_{b,n_i}=\pi^{-1}(T_i).$$
Note that the curves $T_i$ are then Zariski dense in $T$.
Our goal is to show that this last case cannot occur, if the point $b\in B$ satisfies assumption
(\ref{assumption}). This will be done in several steps.

 We claim  first that this  fibration
  $\pi:D\dashrightarrow T$   is preserved by $\phi^l$.  Indeed, the fibration
 is determined by the
  data of the  $2$-form
 $\omega_D$ up to a scalar :  the fibers of $\pi$ are the integral leaves of the
 vector field defined generically on $D$ as the kernel of $\omega_D$. As we have
 $(\phi^l)^*\omega_D=(-2)^l\omega_D$, the claim follows.

   As $D$ is defined only up to bimeromorphic
 transformations, we may assume $\pi$ is actually a morphism. We will denote by
 $\psi:T\dashrightarrow T$ the rational map induced by $\phi^l$ on the basis of the fibration.

 \begin{lemm}\label{utile1} The generic fiber of $\pi$ is a curve of general type.
 \end{lemm}
{\bf Proof.} As the $T_i$ are Zariski dense in $T$, it suffices to prove the same result for the fibration
$\pi:\Sigma'_{b,n_i}\rightarrow T_i$.

Recall now that the surfaces $\Sigma'_{b,j}$ are birationally equivalent to an abelian surface
$A_b$,
and that we proved in section \ref{preli} that we could assume $A_b$ to be
geometrically simple.
We then simply use the fact that curves in  a simple abelian surface are of general type.
\cqfd
We now have:
\begin{lemm}\label{utile2} There exists a rational map
$m:T\dashrightarrow Z$ to a curve $Z$, satisfying the following properties:
\begin{enumerate}
\item \label{item1}
 The restriction of the fibration $\pi:D\rightarrow T$ to a generic fiber $T_z=m^{-1}(z)$
is an isotrivial fibration.
\item \label{item2}The morphism $\phi^l$, which was already proved to descend to $T$, preserves
 the fibration
$m$, acting trivially on the basis $Z$. Thus $\phi^l$ acting on $D$
 acts on each of  the fibers of $m\circ \pi:D\dashrightarrow Z$.
 \end{enumerate}
 \end{lemm}
 {\bf Proof.} Let $g\geq 2$ be the genus of the
 generic fiber of $\pi$. We take for $m$ the rational map to $\mathcal{M}_g$ associated to the fibration
 $\pi$. By definition, the fibration becomes
  isotrivial when restricted to the fibers of $m$, proving \ref{item1}.

 As $\phi^l$ acts on $D$ preserving the fibration $\pi$, and the fibers of $\pi$ are of general type,
 $\phi^l$ induces for generic $t\in T$ an isomorphism
 $$C_t\cong C_{\psi(t)},\,t\in T,$$
 where $C_t:=\pi^{-1}(t)$. This isomorphism shows that
 $m\circ \psi=m$. To prove \ref{item2}, it thus only remains to show that
 $Im\,m$ is a curve.

As the rational self-map $\psi$ is not of finite order on $T$ (because the $T_i$ are of the form
$\psi^k(T_1)$ and are Zariski dense in $T$), the equality
$m\circ \psi=m$ implies that $m$ is not generically finite-to-one.

Thus either $m$ is constant or its image $Z$ has dimension $1$.
In fact $m$ cannot be constant. Indeed, consider the restriction of $m$ to any of the curve $T_i$.
Recall that we have $ \pi^{-1}(T_i)=\Sigma'_{b,n_i}$, which is birationally equivalent to an abelian surface
$A_b$.
If $m$ were constant on $T_i$, then the rational fibration
$A_b\dashrightarrow T_i$ would be isotrivial with fiber $C$ of general type. This does not exist
on any abelian surface $A$, because otherwise, $A$ would be dominated by a product
$C\times C'$, which would send the $C\times c'$
to the fibers of the fibration. But a rational map
$$\alpha:C\times C'\dashrightarrow A$$
is necessarily a morphism which is a ``sum morphism'', that is of the form
$$\alpha(c,c')=\alpha_1(c)+\alpha_2(c').$$
Thus the fibers of this isotrivial fibration would be translates
in $A$ (in the directions given by
$\alpha_2(C')$)  of a given
curve of general type. But two fibers must intersect because
they are of general type, and they can intersect only along the $0$-dimensional
  indeterminacy locus of the considered fibration on $A$. Their intersection
  is thus stable under translation by $\alpha_2(C')$, which is a contradiction.

 \cqfd
 Let $z\in Z$ be a point defined over a number field.
  We will denote by $T_z\subset T$ the corresponding fiber of
 $m$. The surface $\pi^{-1}(T_z)$ admits an isotrivial fibration
 with fiber $C_z$ defined over a number field $k$, and thus is dominated by
 a product $C_z\times_k T_z'$, where $T_z'$ is a finite cover of $T_z$.

For each point $y\in T'_z$ defined over a number field $k(y)$, the curve
$D_y:=\pi^{-1}(y)$ is isomorphic to
$C_z$ and the restricted cycle
$$P'_{y}:=P'_{\mid D_y\times X}\subset D_{y}\times X_{k(y)}\cong C_{z,k(y)}\times X_{k(y)}$$ is
 homogically trivial in the sense that its {\'e}tale cohomology class
 in $H^6_{et}(D_{y,\overline{k}}\times X_{\overline{k}}(3))$ is zero.
 Thus it admits an Abel-Jacobi $l$-adic invariant
 $$\gamma_y\in
 H^1(k(y),H^5(D_{y,\overline{k(y)}}\times X_{\overline{k(y)}}(3))$$
$$ =H^1(k(y),H^5(C_{z,\overline{k(y)}}\times X_{\overline{k(y)}}(3))$$
 which is obtained as in the previous section by considering the continuous {\'e}tale cycle class of
 $P'_y$ in $H^6(D_y\times X_{k(y)},\mathbb{Q}_l(3))$
 and applying the Hochschild-Serre spectral sequence to it.

We choose a point $z\in Z$ defined over a number field
and which is sufficiently general (so as to avoid singular fibers).
 We will prove the following two propositions:
 \begin{prop} \label{prop1fin} The class $\gamma_y$ does not depend on $y\in T'_z(\overline{\mathbb{Q}})$.
 \end{prop}
This makes sense, as for $y,\,y'\in T'_z(\overline{\mathbb{Q}})$, one may choose a common definition field
$k$ for $y$ and $y'$, on which the isomorphism
$D_y\cong D_{y'}\cong C_z$ is also defined. Then
we compare both classes in
$H^1(k,H^5(C_{z,\overline{k}}\times X_{\overline{k}},\mathbb{Q}_l)(3))$.

\vspace{0,5cm}

{\bf Proof of Proposition \ref{prop1fin}.} Indeed, consider the surface
$\pi^{-1}(T_z)$ which admits as a rational finite cover $C_z\times_k T'_z$. Let us denote
by $r:C_z\times_k T'_z\dashrightarrow F$ the natural rational map.
We claim that $(r,Id)^*P'\in CH^3(C_z\times _kT'_z\times X)_\mathbb{Q}$ is a cycle homologous to $0$,
which means that its  {\'e}tale cohomology class on
$C_{z,\overline{k}}\times _{\overline{k}}T'_{z,\overline{k}}
\times _{\overline{k}}X_{\overline{k}}$ or equivalently its Betti cohomology class
on $C_{z,\mathbb{C}}\times T'_{z,\mathbb{C}}
\times X_{\mathbb{C}}$ vanishes.
Indeed, we know that the $2$-form $\omega=P^{'*}\alpha$ on $F$ has the property that
its restriction $\omega_D$ to $D$ is of the form $\pi^*\omega_T$ for some holomorphic
$2$-form on $T$. Of course $\omega_T$ vanishes on the curve $T_z$.
It follows that $\omega_D$ vanishes on the surface $\pi^{-1}(T_z)$.
Let us work in the complex setting: by lemma
\ref{hodge},  the Hodge structure on $H^4(X,\mathbb{Q})_{prim}$
is simple.
As the morphism
$$r^*\circ P^{'*}:H^4(X_\mathbb{C},\mathbb{Q})_{prim}\rightarrow H^2(C_{z,\mathbb{C}}\times
T'_{z,\mathbb{C}},\mathbb{Q})$$
is not injective, because its complexification vanishes on the $(3,1)$-part,
it must be $0$, which proves the claim.

Having this, we conclude that the cycle
$(r,Id)^*P'\in CH^3(C_z\times _kT'_z\times X)$ has an $l$-adic Abel-Jacobi invariant
$$\gamma\in
H^1(k,H^5(C_{z,\overline{k}}\times_{\overline{k}}
 T'_{z,\overline{k}}\times _{\overline{k}} X_{\overline{k}},\mathbb{Q}_l)(3))$$
 which in fact lies in
 $$H^1(k,H^1(C_{z,\overline{k}}\times_{\overline{k}}
 T'_{z,\overline{k}},\mathbb{Q}_l)\otimes _{\mathbb{Q}_l}H^4( X_{\overline{k}},\mathbb{Q}_l)(3))$$
 because $P'$ is obtained from $P$  by applying a Chow-K\"unneth projection on $X$.

We have $$H^1_{et}(C_{z,\overline{k}}\times_{\overline{k}}
 T'_{z,\overline{k}},\mathbb{Q}_l)\cong H^1_{et}(C_{z,\overline{k}},\mathbb{Q}_l)\oplus H^1_{et}(
 T'_{z,\overline{k}},\mathbb{Q}_l),$$
 and thus the first projection of $\gamma$ gives us a class
 $$\gamma_z\in
 H^1(k,H^1(C_{z,\overline{k}},\mathbb{Q}_l)
 \otimes _{\mathbb{Q}_l}H^4( X_{\overline{k}},\mathbb{Q}_l)(3))$$
 which  by definition restricts to
 $\gamma_y$, for any $y\in T'_z$.
\cqfd

\begin{prop} \label{prop2fin}  The class
 $\gamma_y,\,y\in T'_z$ is non zero, under the assumption (\ref{assumption})
made on $b\in B$.
 \end{prop}
{\bf Proof.} The fiber $T_z$ meets the countably many curves $T_i$, as we proved in the proof of lemma
\ref{utile2}.
Let $y\in T_1\cap T_z$ be defined over a number field $k$. The class
$\gamma_y$ is equal to the class $\gamma_{z}$ by Proposition \ref{prop1fin}.
Now, note that, as $y\in T_1$, the fiber $D_{y}$ is a curve in the surface
$\Sigma'_{b,n_1}=\pi^{-1}(T_1)$. This surface is birationally equivalent
to a simple abelian surface, and thus we have an injective
morphism of {\'e}tale cohomology groups:
$$ H^1_{et}(\Sigma'_{b,n_1,\overline{k}},\mathbb{Q}_l)\rightarrow H^1_{et}(D_{y,\overline{k}},
\mathbb{Q}_l),$$
where $k$ is a common field of definition of $T_1$ and $y$, and $\overline{k}$ is
its Galois closure.
This restriction map is $Gal\,(\overline{k}/k)$-equivariant and it makes
the left hand side into a direct summand of right hand side, using the existence
of a non degenerate $Gal\,(\overline{k}/k)$-invariant
intersection pairing on the right hand side.
It follows that there is an induced injection
$$H^1(k, H^1(\Sigma'_{b,n_1,\overline{k}},\mathbb{Q}_l)
\otimes_{\mathbb{Q}_l}H^4(X_{\overline{k}},\mathbb{Q}_l)(3))$$
$$
\hookrightarrow
H^1(k, H^1(D_{y,\overline{k}},\mathbb{Q}_l)
\otimes_{\mathbb{Q}_l}H^4(X_{\overline{k}},\mathbb{Q}_l)(3)),$$
which obviously sends the Abel-Jacobi $l$-adic invariant of the cycle
$P'_{\Sigma'_{b,n_1}}$ to the
$l$-adic invariant of the cycle $P'_{y}$, that is $\gamma_{y}$.
We already proved in proposition \ref{plusfort} that the first one is non zero,
because it
is equal to $(-2)^{n_1}\gamma_{P'_b}$ by property (\ref{crucialpropertyphi}), hence
$\gamma_{y}\not=0$.
\cqfd

Let us now explain how to use these propositions
to get a contradiction, which will exclude the case considered in this section.

The fiber $T_z$ meets the countably many curves $T_i$, as we proved in the proof of Lemma
\ref{utile2}. This fiber is acted on by $\psi$, as proved in lemma \ref{utile2}.
Let $y\in T_1\cap T_z$ be defined over a number field $k$.

By proposition \ref{prop1fin}, its
$l$-adic Abel-Jacobi invariant
$$\gamma_{y}\in H^1(k,H^5(D_{y,\overline{k}}\times X_{\overline{k}})(3))$$
 must be equal
to the $l$-adic Abel-Jacobi invariant
$$\gamma_{y'}\in H^1(k,H^5(D_{y',\overline{k}}\times
_{\overline{k}} X_{\overline{k}})(3))$$
of the cycle
$P'_{y'}\subset D_{y'}\times X$, where $y'=\psi(y)\in \psi(T_1)\cap T_z$,
and $D_{y}$ and $D_{y'}$ are identified via the trivialisation
of the fibration of $\pi$ over $T_z$. Up to replacing the power $l$ by an higher power,
we may assume that this identification is given by
$\phi^l:D_{y}\cong C_{y'}$ because the automorphisms group of these curves is finite.
Thus we get
\begin{eqnarray}\label{autre1}{\phi^l}^*\gamma_{y'}=\gamma_{y},
\end{eqnarray}
where we see $\phi^l$ as a map from $D_{y}$ to $D_{y'}$.

On the other hand, we
 get  from
(\ref{crucialpropertyphi}) and  Bloch-Srinivas decomposition theorem
 (cf \cite{blS} and
\cite{voisinbook}, Corollary 10.20)  the equality
$${\phi^l}^*(P'_{y'})=(-2)^lP'_y+\Gamma\,\,{\rm in}\,\, CH^3(D_{y}\times X)_\mathbb{Q},$$
 where $\Gamma\subset D_{y}\times X$ is a ``vertical'' cycle, supported over
a divisor of $D_{y}$.
As vertical cycles have trivial $l$-adic Abel-Jacobi invariant in the
considered group, this  implies  that
\begin{eqnarray}\label{autre2}{\phi^l}^*\gamma_{y'}=(-2)^l\gamma_{y}.
\end{eqnarray}
As
$$0\not=\gamma_{y}\,\,{\rm in}\,\,
 H^1(k,H^5(D_{y,\overline{k}}\times_{\overline{k}}
 X_{\overline{k}},\mathbb{Q}_l)(3))$$
by proposition \ref{prop2fin},
(\ref{autre1}) and (\ref{autre2}) provide a contradiction.
\cqfd

\end{document}